\documentclass[10pt]{amsart}
\usepackage{amsmath,amssymb}
\usepackage{amsrefs}

\numberwithin{equation}{section}

\newtheorem{theorem}[equation]{Theorem}
\newtheorem{lemma}[equation]{Lemma}

\theoremstyle{definition}
\newtheorem{definition}[equation]{Definition}

\theoremstyle{remark}
\newtheorem{remark}[equation]{Remark}

\begin{document}

\title{Capability of nilpotent products of cyclic groups II}
\author{Arturo Magidin}

\address{Department of Mathematics, University of %
Louisiana--Lafayette, 217 Maxim Doucet Hall, P.O.~Box 41010, Lafayette %
LA 70504-1010} 

\email{magidin@member.ams.org}

\begin{abstract}
In Part~I it was shown that if $G$ is a $p$-group of class~$k$,
generated by elements of orders $1<p^{\alpha_1}\leq \cdots \leq p^{\alpha_r}$,
then a necessary condition for the capability of~$G$ is that $r>1$ and
$\alpha_r\leq \alpha_{r-1}+\lfloor\frac{k-1}{p-1}\rfloor$. It was also shown
that when $G$ is the $k$-nilpotent product of the cyclic groups
generated by those elements and $k=p=2$ or $k<p$, then
the given conditions are also sufficient.
We make a correction related to the small class
case, and extend the sufficiency result to $k=p$ for arbitrary prime~$p$.
\end{abstract}

\maketitle

Recall that a group $G$ is said to be \textit{capable} if and only if
there exists a group~$H$ such that $G\cong H/Z(H)$, where $Z(H)$ is
the center of~$H$. In~\cite{capable} we proved that if $G$ is a
capable $p$-group of class~$k$, generated by $x_1,\ldots,x_r$, with
$x_i$ of order $p^{\alpha_i}$, $1\leq \alpha_1\leq\cdots\leq \alpha_r$, then $r>
1$ and $\alpha_r\leq \alpha_{r-1}+\lfloor\frac{k-1}{p-1}\rfloor$.  We also
proved that if~$G$ is the $k$-nilpotent product of the cyclic $p$-groups
generated by the $x_i$,  then the conditions are also
sufficient for the cases $k<p$ and $k=p=2$.

The purpose of this note is twofold: first, we will note an error in a
lemma that was used in the proof of the small class case and make the
necessary corrections to justify that result. Second, we will extend
the result to the case $k=p$ with $p$ an arbitrary prime.  Since we
follow closely on~\cite{capable}, we refer the reader there for the
relevant definitions and conventions. 

I am extremely grateful to Prof.~T. C.~Hurley who brought to my attention
the results from \cites{basiccomms,hurley}; these results
allowed the correction of the error noted above, as well as
simplifying my argument for the $k=p$ case.

\section{Shoving commutators}

In \cite{capable}, the last clause of Lemma~4.2(ii) is incorrect.
Because of this error, the last assertion in Lemma~4.3 is also
incorrect; the proof of Theorem~4.4, which describes the center of a
$k$-nilpotent product of cyclic $p$-groups when $k\leq p$, relied on
that incorrect assertion and so has a gap. In this section we will
provide the necessary correction to justify the conclusion of that
theorem. Once it is established, the rest of the proof of the
small class case will follow.

The error in question is the following: we start with the free group
$F$ on $x_1,\ldots,x_r$, and a basic commutator $[u,v]$ of weight
equal to $k\geq 2$. Then we considered $[u,v,x_r]$; when $v\leq x_r$,
this is a basic commutator. If $v>x_r$, then we rewrite $[u,v,x_r]$
modulo $F_{k+2}$ as $[u,x_r,v][v,x_r,u]^{-1}$. The incorrect clause
asserted that this expresses $[u,v,x_r]$ modulo~$F_{k+2}$ as a
product of basic commutators and their inverses, but this is not
necessarily the case; there is no warrant for asserting that $[u,x_r]$
or $[v,x_r]$ will necessarily be basic commutators (though they are
for small values of~$k$), nor that $[v,x_r]>u$, another
requirement. The main idea is sound, however, if one continues this
process until the resulting expression consists only of powers of
basic commutators. Fortunately, while the final expression may be
too complex to describe in general, one does have control over the
smallest commutator that occurs in that expression, as was shown in
\cites{basiccomms,hurley}.

\begin{definition} Let $c$ be a basic commutator on
  $x_1,\ldots,x_r$. Then ${\rm wt}(c)$ is the weight of $c$ in
  the~$x_i$.
\end{definition}

\begin{definition}[\cite{basiccomms}*{\S 18}]
 Let $u$ and $v$ be basic commutators on
  $x_1,\ldots,x_r$. The commutator $[u{\leftarrow} v]$ is defined
  recursively as:
\begin{itemize}
\item[(i)] If $v>u$, then $[u{\leftarrow} v] = [v{\leftarrow} u]$.
\item[(ii)] If $v<u$ and $u=[c_1,c_2]$ and $c_2>v$, then $[u{\leftarrow}
  v] = [c_1{\leftarrow} v,c_2]$.
\item[(iii)] Otherwise, $[u{\leftarrow} v] = [u,v]$.
\end{itemize}
\end{definition}

\noindent Here is an explicit description, easily established:

\begin{lemma}[cf. \cite{witts}*{Lemma~2.4}]
Let $u$ and $v$ be basic commutators on $x_1,\ldots,x_r$ with $u>v$.
If ${\rm wt}(u)=1$, then $[u{\leftarrow} v]=[u,v]$. If ${\rm
  wt}(u)>1$, then letting $u = [c_1,\ldots,c_n]$ where $c_i$ is a basic
commutator and ${\rm wt}(c_1)=1$, we have:
\begin{itemize}
\item[(1)] If $c_2>v$, then $[u{\leftarrow} v] =
  [c_1,v,c_2,\ldots,c_n]$.
\item[(2)] Otherwise, $[u{\leftarrow} v] =
  [c_1,\ldots,c_j,v,c_{j+1},\ldots,c_n]$, where $j$ is the largest
  index such that $c_j\leq v$.
\end{itemize}
\label{lemma:description}
\end{lemma}

Given a commutator $c=[r,s]$, we will say informally that $r$ is
the ``left entry'' of~$c$, and that $s$ is the ``right entry'' of~$c$.

The following lemma is due to Ward (modulo a different 
definition of the order among basic commutators of the same
weight); the proofs are straightforward.

\begin{lemma}[\cite{basiccomms}*{Lemma~18.1}]
Let $u$ and $v$ be basic commutators on $x_1,\ldots,x_r$, and assume
that $u>v$.
\begin{itemize}
\item[(i)] $[u{\leftarrow} v]$ exists and is a basic commutator on
  $x_1,\ldots,x_r$.
\item[(ii)] ${\rm wt}([u{\leftarrow} v]) = {\rm wt}(u)+{\rm wt}(v)$.
\item[(iii)] $u < [u{\leftarrow} v]$.
\item[(iv)] If ${\rm wt}(u)>1$, then the right entry of $[u{\leftarrow} v]$ is equal to the
  larger of $v$ and the right entry of~$u$.
\item[(v)] If $v<u_1<u_2$, then $[u_1{\leftarrow}
  v]<[u_2{\leftarrow} v]$.
\item[(vi)] If $v_1<v_2<u$, then $[u{\leftarrow}
  v_1]<[u{\leftarrow} v_2]$.
\item[(vii)] If $[u_1{\leftarrow} v] = [u_2{\leftarrow} v]$, then
  $u_1=u_2$.
\end{itemize}
\label{lemma:propshove}
\end{lemma}

We can think of $[u{\leftarrow} v]$ as the basic commutator which
results from ``shoving'' $v$ into its correct position inside of~$u$
(hence the title of this section).

Some of the results above, in particular (v) and (vi), were also
obtained independently by Waldinger~\cite{naturalorder}. 

The following result is also essentially contained in
\cites{basiccomms,naturalorder}. However, both authors use a
prefered ordering among basic commutators that is different from ours,
so their conclusions also read differently. Because of this, we
provide a proof.

\begin{lemma}[cf.~\cite{basiccomms}*{Lemma 18.3}, \cite{hurley}*{Lemma
      1.2}] Let $F$ be the free group on $x_1,\ldots,x_r$. Let $u$ and
$v$ be basic commutators and let $k={\rm wt}(u)+{\rm
wt}(v)$. Then
\[ [u,v] \equiv [u{\leftarrow} v]^{\epsilon}
c_1^{\alpha_1}c_2^{\alpha_2}\cdots c_m^{\alpha_m} \pmod{F_{k+1}},\]
where $F_{k+1}$ is the $(k+1)$st term of the lower central series of
$F$, $\epsilon=\pm1$, $\alpha_i$ are integers, each $c_i$ is a
basic commutator of weight~$k$, and $[u{\leftarrow}
  v]<c_1<\cdots<c_m$. Moreover, if $u\geq v$, then we may choose $\epsilon=1$.
\label{lemma:basiclemma}
\end{lemma}


\begin{proof} It is enough to establish the result when $u>v$: if $u=v$, then
$[u,v]$ and $[u{\leftarrow} v]$ are both trivial, so setting $\epsilon=1$ and $m=0$ proves
the result. And if
$v>u$, then $[u,v] = [v,u]^{-1}$; assuming the result holds when the
left entry is greater than the right entry, and since $F_k/F_{k+1}$ is
abelian, we obtain
\begin{alignat*}{2}
[u,v] = [v,u]^{-1} &\mathop{\equiv} \left([v{\leftarrow}
  u]c_1^{\alpha_1}\cdots c_n^{\alpha_n}\right)^{-1}
  &\pmod{F_{k+1}}\hphantom{.}\\
& \mathop{\equiv} \left[v{\leftarrow} u\right]^{-1}c_1^{-\alpha_1}\cdots
  c_n^{-\alpha_n} &\pmod{F_{k+1}}\hphantom{.}\\
& \mathop{\equiv} \left[u{\leftarrow} v\right]^{-1}c_1^{-\alpha_1}\cdots
  c_n^{-\alpha_n}&\pmod{F_{k+1}}.
\end{alignat*}
So we assume without loss of generality that $u>v$.

We proceed by induction on $k$. If $k=2$, then $u=x_j$ and $v=x_i$
with $i < j$; hence $[u,v]=[u{\leftarrow} v]$. If $k=3$, then
$u=[x_j,x_i]$ and $v=x_{\ell}$, with $1\leq i<j\leq r$; if $i\leq \ell$ then
$[u,v] = [u{\leftarrow} v]$ and we are done. If $i>\ell$, then from
\cite{capable}*{Prop.~2.2(iv)} we have that $[u,v] \equiv
[x_j,x_{\ell},x_i][x_i,x_{\ell},x_j]^{-1}\pmod{F_4}$, and $[u{\leftarrow}
v]=[x_j,x_{\ell},x_i]$, which is strictly smaller than $[x_i,x_{\ell},x_j]$, so
the result also holds.

Assume that $k>3$ and the result is true for all commutators $[c_1,c_2]$
where $c_1$ and~$c_2$ are basic commutators with ${\rm wt}(c_1)+{\rm
wt}(c_2)<k$ and $c_1>c_2$. We will now argue by ``descending induction''
on~$v$. Picking the largest possible weight for $v$ for which ${\rm
wt}(u)+{\rm wt}(v)=k$ and $u>v$ yields ${\rm wt}(u)={\rm wt}(v)$ or
${\rm wt}(u)={\rm wt}(v)+1$. Write $u=[a,b]$ (which we can do since
${\rm wt}(u)\geq 2$).  Then ${\rm wt}(b)\leq \frac{1}{2}{\rm
wt}(u)$. If ${\rm wt}(b)\geq{\rm wt}(v)$, then we would have ${\rm
wt}(v)\leq \frac{1}{2}{\rm wt}(u)\leq \frac{1}{2}({\rm wt}(v)+1)$; for
this to hold we must have ${\rm wt}(v)=1$ and ${\rm wt}(u)\leq 2$,
contradicting the assumption that $k>3$.  Hence ${\rm wt}(b)<{\rm
wt}(v)$, so $b<v$ and $[u,v] = [u{\leftarrow} v]$; thus the result holds
in this case.

Suppose then the result holds for $[c_1,c_2]$ whenever $c_1$ and~$c_2$ are
  basic commutators, $c_1>c_2$, and either ${\rm wt}(c_1)+{\rm
  wt}(c_2)<k$ or ${\rm wt}(c_1)+{\rm wt}(c_2)=k$ and $c_1>c_2>v$.
Write $u=[a,b]$; if $b\leq v$, then $[u,v]=[u{\leftarrow} v]$ and
we are done. Otherwise, again from 
\cite{capable}*{Prop.~2.2(iv)} we have that
$[u,v] \equiv [a,v,b][b,v,a]^{-1} \pmod{F_{k+1}}$.
Let $\kappa = {\rm wt}(a)+{\rm wt}(v)$ and $\lambda={\rm
  wt}(b)+{\rm wt}(v)$. By induction, we know that:
\[
\relax
[a,v] \equiv [a{\leftarrow} v]\prod_{i=1}^r
c_i^{\alpha_i}\pmod{F_{\kappa+1}},\quad\mbox{and}\quad
[b,v]  \equiv  [b{\leftarrow} v]\prod_{j=1}^s
d_j^{\beta_j}\pmod{F_{\lambda+1}},
\]
where $\alpha_i,\beta_j$ are integers, the $c_i$ are basic commutators of
 weight $\kappa$, the $d_j$ are basic commutators of weight $\lambda$, and
 the inequalities $[a{\leftarrow} v]<c_1<\cdots<c_r$ and $[b{\leftarrow}
 v]<d_1<\cdots<d_s$ hold.  Since $\kappa+{\rm wt}(b) = \lambda+{\rm
 wt}(a)=k$, from well-known commutator identities (e.g., those in
 \cite{capable}*{Prop 2.2}) we obtain:
\[ [u,v]  \equiv  [a{\leftarrow} v,b]\left(\prod_{i=1}^r
[c_i,b]^{\alpha_i}\right)[b{\leftarrow} v,a]^{-1}\left(\prod_{j=1}^s
 [d_j,a]^{-\beta_j}\right)\pmod{F_{k+1}}.\] 
Note that $[a{\leftarrow}
 v,b] = [u{\leftarrow} v]$, that ${\rm wt}(c_i)+{\rm wt}(b)={\rm
 wt}(d_j)+{\rm wt}(a)=k$, and likewise ${\rm wt}([b\leftarrow
 v])+{\rm wt}(a)=k$.
Our
 result will therefore follow if we can prove that each of $[c_i,b]$,
 $[b{\leftarrow} v,a]$ and $[d_j,a]$ is congruent modulo $F_{k+1}$
 to a product of powers of basic commutators of weight $k$, 
 each strictly larger than $[u{\leftarrow} v]$; then we can invoke the
 fact that $F_k/F_{k+1}$ is abelian to obtain an expression for
 $[u,v]$ modulo~$F_{k+1}$ of the desired form. We remove any commutators that are
 trivial, and consider each of the remaining ones in turn.

Since $b>v$, the induction hypothesis allows us
to rewrite each $[c_i,b]$ as a product of powers of basic commutators,
each greater than or equal to $[c_i{\leftarrow} b]$. We know that
$c_i>[a{\leftarrow} v]$, so we have that $[c_i{\leftarrow} b] > [
[a{\leftarrow} v]{\leftarrow} b]$; since $b$ is no smaller than the right
entry of $[a{\leftarrow} v]$ we know that $[[a{\leftarrow} v]{\leftarrow} b]
= [a{\leftarrow} v,b]=[u{\leftarrow} v]$. Therefore, $[c_i{\leftarrow}
b]>[u{\leftarrow} v]$ and so all basic commutators that appear in the
expression for $[c_i,b]$ are also strictly larger than $[u{\leftarrow}
v]$. So we can certainly replace each of the $[c_i,b]$ as
needed. 
  
If $[b{\leftarrow} v]<a$, then we replace $[b{\leftarrow} v,a]^{-1}$ with
$[a,b{\leftarrow} v]$. Since $b<[b{\leftarrow} v]$ and the right
entry of $a$ is less than or equal to $b$, we have that
$[a,b{\leftarrow} v]$ is a basic commutator; also since
$b<[b{\leftarrow} v]$ we deduce that $[u{\leftarrow} v]=[a{\leftarrow}
v,b]<[a,b{\leftarrow} v]$. On the other hand, if $[b{\leftarrow}
v]>a$, then we know that the right entry of $[b{\leftarrow} v]$ is
strictly less than $b$ (equal to either $v$ or to the right entry
of~$b$), hence strictly smaller than~$a$; thus, $[b{\leftarrow} v,a]$
is already a basic commutator. The right entry of this latter
commutator is $a$, so $[b{\leftarrow} v,a]>[a{\leftarrow}
v,b]=[u{\leftarrow} v]$. This shows the commutator $[b{\leftarrow}
  v,a]$ is either a basic commutator greater than $[u{\leftarrow} v]$,
or the inverse of a basic commutator greater than $[u{\leftarrow}
  v]$. 

Finally, we come to the commutators $[d_j,a]$. If $a>d_j$, then we
replace $[d_j,a]$ with $[a,d_j]^{-1}$. Since $d_j>[b{\leftarrow}
v]>b$, it follows that the right entry of $a$ is strictly smaller
than $d_j$, so $[a,d_j]$ is a basic commutator; and $d_j>b$ also
implies that $[a,d_j]>[a{\leftarrow} v,b]=[u{\leftarrow} v]$. On the
other hand, if $d_j>a$, since $a>v$ we can again apply induction to replace
$[d_j,a]$ with a product of $[d_j{\leftarrow} a]$ 
times powers of basic commutators strictly larger than $[d_j{\leftarrow}
a]$. The right entry of $[d_j{\leftarrow} a]$ is no smaller than
$a$, and hence is strictly larger than $b$, the right entry of
$[u{\leftarrow} v]$. Thus, we can also replace each $[d_j,a]$ with a product of
powers of basic commutators, each larger than $[u{\leftarrow} v]$. This
proves the lemma. 
\end{proof}

\begin{lemma}[cf.~\cite{hurley}*{Lemma~1.3}]
Let $F$ be the absolutely free group on $x_1,\ldots,x_m$, and suppose
that $c\equiv b_1^{\alpha_1}\cdots b_t^{\alpha_t}\pmod{F_{k+1}}$,
where $\alpha_i$ are integers, $b_i$ are basic commutators of weight
exactly $k$, and $b_1<b_2<\cdots<b_t$. If $a$ is a basic commutator of
weight $\ell$,  then
$[c,a] \equiv [b_1{\leftarrow} a]^{\pm\alpha_1} u_1^{\beta_1}\cdots
u_s^{\beta_s}\pmod{F_{k+\ell+1}}$ where the $\beta_i$ are integers,
$u_i$ are basic commutators with ${\rm wt}(u_i)=k+\ell$, and 
$[b_1{\leftarrow} a]<u_1<\cdots<u_s$. Moreover, if $c>a$
then the exponent of $[b_1{\leftarrow} a]$ may be taken to be $\alpha_1$.
\label{lemma:generalshoving}
\end{lemma}

\begin{proof} We have that $[c,a]\equiv
  [b_1,a]^{\alpha_1}\cdots[b_t,a]^{\alpha_t}\pmod{F_{k+\ell+1}}$;
  since the $b_i$ are in increasing order, the corresponding
  $[b_i{\leftarrow} a]$ are also in increasing order; the result now
  follows from the fact that $F_{k+\ell}/F_{k+\ell+1}$ is abelian and from
Lemma~\ref{lemma:basiclemma}.
\end{proof}

With this result, we can replace the argument based on the erroneous
Lemma~4.3 and prove Theorem~4.4 from \cite{capable}:

\begin{theorem}[\cite{capable}*{Theorem~4.4}]
For a positive integer $k$ and a prime~$p$ with $p\geq k$, let
$C_1,\ldots,C_r$ be cyclic $p$-groups generated by $x_1,\ldots,x_r$
respectively, with $p^{\alpha_i}$ being the order of~$x_i$, and assume
that $1\leq \alpha_1\leq\cdots\leq \alpha_r$. If $G$ is the
$k$-nilpotent product of the $C_i$, $G=C_1\amalg^{\germ
  N_k}\cdots\amalg^{\germ N_k} C_r$, then $Z(G)=\langle
x_r^{p^{\alpha_{r-1}}},G_k\rangle$. 
\label{thm:centerknil}
\end{theorem}

\begin{proof} The center contains both $x_r^{p^{\alpha_{r-1}}}$ and
  $G_k$ by \cite{capable}*{Lemma 3.11} and the properties of the nilpotent
  product. The reverse inclusion is established by induction, the case
  $k=1$ being trivial and the case $k=2$ having been proven in
  \cite{capable}*{Lemma~4.1}. If we consider $G/G_k$ we obtain the
  $(k-1)$-nilpotent product of the $C_i$, from which we have that
  $\langle x_r^{p^{\alpha_{r-1}}},G_k\rangle\subseteq Z(G)\subseteq
  \langle x_r^{p^{\alpha_{r-1}}}, G_{k-1}\rangle$ by induction. To
  prove equality, it is enough to show that if $g\in G_{k-1}\cap
  Z(G)$, then $g\equiv e\pmod{G_k}$. Write $g\equiv
  b_1^{\beta_1}\cdots b_t^{\beta_t}\pmod{G_k}$, with each $b_i$ a
  basic commutators of weight $k-1$, and $b_1<\cdots<b_t$; from
  \cite{struikone}*{Theorem~3} we know this expression is unique if we
  require that each $\beta_i$ satisfy $0\leq \beta_i <
  p^{\alpha_{s_i}}$, where $s_i$ is the smallest index of a generator
  that appears in the full expression for $b_i$. If $t=0$ then
  trivially $g\in G_k$. Assuming the conclusion holds for expressions
  with fewer terms, by
  Lemma~\ref{lemma:generalshoving} we have that $e=[g,x_r] =
  [b_1{\leftarrow} x_r]^{\beta_1}\prod u_j^{\gamma_j}$ (equality since
  $G_{k+1}$ is trivial), where $u_j$ are basic commutators of
  weight~$k$, with $[b_1\leftarrow x_r]<u_1<u_2<\cdots$. Again, by the
  normal form proven in \cite{struikone}*{Theorem~3}, and since the
  order of $[b_1{\leftarrow} x_r]$ must be equal to the order of
  $b_1$, we deduce $\beta_1=0$ so we may express $g$ modulo $G_k$
  using fewer than~$t$ powers of basic commutators, and by induction
  we deduce $g\in G_k$, as claimed.
\end{proof}

\section{The case $k=p$}

In this section we will extend the main result from \cite{capable}
to the case $k=p$ with $p$ an arbitrary prime. We will do so by
computing the center of a $(p+1)$-nilpotent product of cyclic
$p$-groups much in the same way as above, using a normal form for the
elements of such a product that was obtained by
R.R.~Struik in her detailed
study~\cite{struiktwo}. Lemma~\ref{lemma:generalshoving} will also play a
key part.

\begin{definition} Let $G$ be a group, and let $x,y\in G$. We define
  $[x,{}_{1}\,y] = [x,y]$ and $[x,{}_{n+1}\,y]=[x,{}_{n}\,y,y]$, where
  $n>1$ is an integer.
\end{definition}

The main difficulty in a straightforward extension of the result lies
in the fact that the basic commutators are no longer a good choice for
a ``basis'' for the normal form in the case of the $(p+1)$-nilpotent
product of cyclic $p$-groups, because there are nontrivial relations
between them; for example, a sufficiently high power of $[b,a]$ will
be nontrivial and equal to a power of $[b,{}_p\, a]$. In order to
bypass this difficulty, one chooses a slightly different set of
distinguished commutators, by replacing the basic commutators
$[b,{}_p\, a]$ and $[b,a,{}_{p-1}\,b]$ with the (nonbasic) commutators
$[b,a^p]$ and $[b^p,a]$, respectively.  The normal form result appears
in \cite{struiktwo}*{Theorem~6}, and is as follows: every element $g$
of the $(p+1)$-nilpotent product of cyclic groups generated by
elements $x_1,\ldots,x_r$, with $x_i$ of order $p^{\alpha_i}$, $1\leq
\alpha_1\leq\cdots\leq \alpha_r$, can be written uniquely as $g=\prod
c_i^{\beta_i}$, where $c_1<c_2<\cdots$ is the sequence of basic
commutators of weight at most $p+1$ in $x_1,\ldots,x_r$, except that
the basic commutator $[x_j,{}_p\,x_i]$ is replaced by the commutator
$v_{ji}'=[x_j,x_i^p]$, and the basic commutator $[x_j,x_i,{}_{p-1}\,
x_j]$ is replaced by the commutator $v_{ji}''=[x_j^p,x_i]$; the
$\beta_i$ are nonnegative integers satisfying $0\leq \beta_i < N_i$,
where:
\begin{equation}
N_i = \left\{\begin{array}{ll}
p^{\alpha_i}&\mbox{if ${\rm wt}(c_i) = 1$ and $c_i=x_i$;}\\
p^{\alpha_k+1}&\mbox{if $c_i = [x_j,x_k]$, $1\leq k<j\leq r$;}\\
p^{\alpha_k-1}&\mbox{if $c_i=v_{jk}'= [x_j,x_k^p]$, $1\leq k<j\leq r$;}\\
p^{\alpha_k-1}&\mbox{if $c_i = v_{jk}'' = [x_j^p,x_k]$, $1\leq k<j\leq
  r$ and $\alpha_k=\alpha_j$;}\\
p^{\alpha_k}&\mbox{if $c_i=v_{jk}'' = [x_j^p,x_k]$, $1\leq k<j\leq r$ and
  $\alpha_k<\alpha_j$;}\\
&\mbox{if $c_i$ is any other basic commutator}\\
p^{\alpha_{s_i}}&\mbox{\ and
  $s_i$ is the smallest index occurring}\\
&\mbox{\ in the full expression for
  $c_i$.}
\end{array}\right.
\label{eq:defofNi}
\end{equation}

\begin{remark} There is a slight inconsistency between the above and
  the statement of \cite{struiktwo}*{Theorem~6}; in the latter, the
  range for the exponents of $[x_j,x_k]$ is not explicitly specified,
  and would be $0$ to $p^{\alpha_k}$ following the general
  case. However, the discussion leading up to the theorem, and in
  particular \cite{struiktwo}*{Equation~60} states that the exponent
  will be taken modulo $p^{\alpha_k+1}$; and this is explicitly the
  case in \cite{struikone}*{Theorem 4} which deals with $p=2$. So it
  seems clear that this is an inadvertent omission in the statement
  of Theorem~6. Nonetheless, our argument will avoid consideration of
  the specific exponent of these commutators except in the case $p=2$.
\end{remark}

We want to describe the center of a $(p+1)$-nilpotent product of
cyclic $p$-groups. The idea is the same one as was used above: if we
let $G$ be the $(p+1)$-nilpotent product of cyclic $p$-groups, then it
is easy to show that $Z(G)$ has upper and lower bounds determined by a
power of $x_r$ and $G_{p+1}$ below, and a power of $x_r$ and~$G_p$
above. At this point we have two extra difficulties we did
not encounter above: the first is that the power of $x_r$ is not the
same in the two bounds, whereas it was the same in the proof of
Theorem~\ref{thm:centerknil}. This can be dealt with in a
straighforward way and we do so first; we will return to the second
difficulty after this lemma:

\begin{lemma}
Let $p$ be a prime, and let $\alpha,\beta$ be positive integers with
$\alpha<\beta$. Let $G=\langle x\rangle \amalg^{\germ N_{p+1}}
\langle y\rangle$, where $x$ generates a cyclic group of order
$p^{\alpha}$ and $y$ generates a cyclic group of order
$p^{\beta}$. Then $[y^{p^{\alpha}},x] = [y^p,x]^{p^{\alpha-1}}$. In
particular, $y^{p^{\alpha}}$ is not central in~$G$.
\label{lemma:notcentral}
\end{lemma}

\begin{proof} All basic commutators in $x$ and~$y$ are of exponent
  $p^{\alpha}$ in $G$, except for $y$ and perhaps $[y,x]$. This can be
  easily established using for example
  \cite{struikone}*{Lemma~H2}. Thus, from \cite{struiktwo}*{Lemma~4}
  we obtain that:
\[ [y^{p^{\alpha}},x] =
  [y,x]^{p^{\alpha}}[y,x,{}_{p-1}\,y]^{\binom{p^{\alpha}}{p}}.\]
Since $[y,x,{}_{p-1}\,y]$ is of exponent $p^{\alpha}$ and
$\binom{p^{\alpha}}{p}\equiv p^{\alpha-1}\pmod{p^{\alpha}}$, we obtain
\[ [y^{p^{\alpha}},x] =
   [y,x]^{p^{\alpha}}[y,x,{}_{p-1}\,y]^{p^{\alpha-1}}.\]
On the other hand, from \cite{struiktwo}*{Equation~(58)} we
have:
\begin{equation}
[y^p,x] = [y,x]^p \left(\prod u_i^{pg_i}\right) [y,x,{}_{p-1}\,y],
\label{eq:ypcommx}
\end{equation}
where the $g_i$ are integers, and $u_i$ are basic commutators of
weight at least three and at most $p+1$ in $x$ and~$y$, omitting both
$[y,x,{}_{p-1}\,y]$ and $[y,{}_{p}\,x]$. From this, since all $u_i$
are of exponent $p^{\alpha}$, we obtain by
\cite{struikone}*{Theorem~H3} that:
\begin{equation}
\relax[y^p,x]^{p^{\alpha-1}} =
   [y,x]^{p^{\alpha}}[y,x,{}_{p-1}\,y]^{p^{\alpha-1}}.
\label{eq:ytothepxraised}
\end{equation}
Therefore, $[y^{p^{\alpha}},x]=[y^p,x]^{p^{\alpha-1}}$, as claimed.
Since we are assuming $\alpha\neq \beta$, the normal form described
   above ensures that $[y^p,x]^{p^{\alpha-1}}\neq e$, and so
   $[y^{p^{\alpha}},x]\neq e$, as claimed.
\end{proof}

The second difficulty alluded to above is a bit more subtle. Once
again the result will come down to proving that if $g\in G_p\cap
Z(G)$, then $g\equiv e \pmod{G_{p+1}}$. If we write $g$ modulo $G_{p+1}$ as a product
of basic commutators of weight exactly $p$ (which can be done since
$G/G_{p+1}$ is the $p$-nilpotent product and the usual normal form
works), and apply
Lemma~\ref{lemma:generalshoving} to compute $[g,x_r]$, we will obtain
$[g,x_r]$ as a product of powers of basic commutators of weight exactly
$p+1$. However, this may not be in the normal form for elements
of~$G$; e.g., 
if any of the basic commutators $[x_j,x_i,{}_{p-1}\, x_j]$ or
$[x_j,{}_{p}\,x_i]$ occur in that expansion then we must 
replace them by expressions using identity (\ref{eq:ypcommx}) and
a similar identity for $[y,x^p]$ \cite{struiktwo}*{Equation~(57)}.
After replacing the occurrences, we must again apply the
collection process to the resulting expression before it will be in
normal form. 

During all of these modifications it might be, at least in principle,
that we modify the exponent of the leading factor in the expression for
$[g,x_r]$ (or even completely replace this leading factor if it is one of the
troublesome basic commutators); thus the argument becomes more
involved. In addition, it may be that the range for the exponents
for the leading factors of $g$ and of $[g,x_r]$ are different.
However, by being careful about just what modifications may
be needed and what they would entail, and sometimes considering
$[g,x_{r-1}]$ instead of $[g,x_r]$,  we can nonetheless push the argument
through to a happy conclusion.

\begin{theorem} Let $p$ be a prime and let $C_1,\ldots,C_r$ be cyclic $p$-groups generated by
$x_1,\ldots,x_r$, of orders
$1<p^{\alpha_1}\leq\cdots\leq p^{\alpha_r}$ 
respectively. If $G$ is the
$(p+1)$-nilpotent product of the $C_i$, $G=C_1\amalg^{\germ
  N_{p+1}}\cdots\amalg^{\germ N_{p+1}} C_r$, then $Z(G)=\langle
x_r^{p^{\alpha_{r-1}+1}},G_{p+1}\rangle$. 
\label{thm:centerpplusone}
\end{theorem}

\begin{proof} That $x_r^{p^{\alpha_{r-1}+1}}$ lies in the center
  follows from \cite{capable}*{Lemma 3.11}; the fact that $G$ is of
  class $p+1$ guarantees that $G_{p+1}\subseteq Z(G)$.

To prove the reverse inclusion, consider $G/G_{p+1}$. By
Theorem~\ref{thm:centerknil} we know the center is generated by (the
images of) $x_r^{p^{\alpha_{r-1}}}$ and $G_{p}$. Pulling back to $G$
we obtain the inclusions
$\langle x_r^{p^{\alpha_{r-1}+1}},G_{p+1}\rangle \subseteq
Z(G)  \subseteq \langle
x_r^{p^{\alpha_{r-1}}},G_{p}\rangle$.
By Lemma~\ref{lemma:notcentral}, if $\alpha_{r-1}<\alpha_r$, then
$x_r^{p^{\alpha_{r-1}}}$ is not central; if $\alpha_{r-1}=\alpha_r$,
then both $x_r^{p^{\alpha_{r-1}}}$ and $x_r^{p^{\alpha_{r-1}+1}}$ are
trivial. So in either case we have
$\langle x_r^{p^{\alpha_{r-1}+1}},G_{p+1}\rangle \subseteq
Z(G) \subseteq \langle
x_r^{p^{\alpha_{r-1}+1}},G_{p}\rangle$.
The theorem will be
established if we can show that for any $g\in G_p$, if $g\in Z(G)$
then $g\equiv e \pmod{G_{p+1}}$. Indeed, if $g\in G_p\cap Z(G)$, then
we can write 
\begin{equation}
g\equiv c_1^{\beta_1}\cdots c_m^{\beta_m} \pmod{G_{p+1}},
\label{eq:expressionforg}
\end{equation}
where $c_1<\cdots<c_m$ are basic commutators of weight exactly $p$,
and $\beta_i$ are nonnegative integer that satisfy $0\leq \beta_i <
p^{\alpha_{s_i}}$, where $s_i$ is the smallest index of a generator
that occurs in the full expression of~$c_i$.
We wish to show
that $g\equiv e\pmod{G_{p+1}}$, and we will do so by induction on~$m$. The result is
trivial if $m=0$; assume then the result holds for all $g$ expressed
as a product of $k$ powers of basic commutators of weight exactly~$p$,
with $0\leq
k<m$. We consider several cases
depending on the nature of the basic commutator $c_1$. 

\textsc{Case 1:} \textit{The right entry of $c_1$ is of weight at
  least two.} Consider $[g,x_r]$. By Lemma~\ref{lemma:generalshoving},
  we have:
\begin{equation}
\relax[g,x_r]=[c_1{\leftarrow} x_r]^{\beta_1}d_1^{\gamma_1}\cdots
  d_n^{\gamma_n},
\label{eq:gcommxr}
\end{equation}
where $\gamma_i$ are integers, the $d_i$ are basic commutators, and
$[c_1{\leftarrow} x_r]<d_1<\cdots<d_n$. We may assume that $0<\gamma_i
< p_{\alpha_{s_i}}$ where $s_i$ is the smallest index of a generator
that occurs in the full expression for $d_i$; this equals the
corresponding $N_i$ from~(\ref{eq:defofNi}) except in the case where
$d_i$ is one of the troublesome commutators.
Since the right entry of $c_1$ is of
weight at least two, so is the right entry of $[c_1{\leftarrow} x_r]$,
and the same holds for each $d_i$. Thus, this expression is already in
normal form and no replacements need to be made. The range of
exponents for $[c_1{\leftarrow} x_r]$ goes from $0$ to
$p^{\alpha_{s}}$, where $s$ is the smallest index that occurs in the
full expression of $[c_1{\leftarrow} x_r]$, which is the same as the
smallest index that occurs in the full expression for $c_1$, namely
$s_1$. Since $g$ is central, we must have $\beta_1\equiv 0
\pmod{p^{\alpha_{s_1}}}$; and from our assumption that $0\leq \beta_1
< p^{\alpha_{s_1}}$ we deduce that $\beta_1=0$. Thus we have $g\equiv
c_2^{\beta_2}\cdots c_m^{\beta_m}\pmod{G_{p+1}}$, and by induction we
deduce that $g\equiv e \pmod{G_{p+1}}$, as desired.

\textsc{Case 2:} \textit{The right entry of $c_1$ is of weight $1$,
and $c_1$ involves at least two generators other than $x_r$.} 
We again consider $[g,x_r]$.
Note that since the right entry of $c_1$ is of weight~$1$, then
  $[c_1{\leftarrow} x_r]=[c_1,x_r]$. Since $[c_1,x_r]<d_i$ for each
  $d_i$ in~(\ref{eq:gcommxr}), the only basic commutators that may
  need to be replaced occur among the $d_i$ and are of the form
  $[x_r,x_i,{}_{p-1}\,x_r]$, which are replaced using 
(\ref{eq:ypcommx}); each of the commutators that are introduced 
  involve only two generators, and so will not equal $c_1$. 
After doing the replacement we must apply the collection process to
rewrite the entire expression in normal form. During the collection,
since in the expression all factors are commutators of weight at least
two, we will only introduce commutators $[b,a]$ in which $a$ is of
weight at least two; again, they will not equal $c_1$. Thus, after
rewriting~(\ref{eq:gcommxr}) in normal form, the exponent of
$[c_1{\leftarrow} x_r]$ will remain $\beta_1$. Since $[g,x_r]=e$, we
must have $\beta_1\equiv 0 \pmod{p^{s_1}}$, which as above yields the
conclusion that $g\in G_{p+1}$, as desired.

\textsc{Case 3:} \textit{The right entry of $c_1$ is of weight $1$,
  and $c_1$ involves only the generators $x_r$ and $x_i$ for some
  $i<r-1$.} Note that we will have $0\leq \beta_1<p^{\alpha_i}$. This
  time we consider $[g,x_{r-1}]$. We have
\begin{equation}
\relax[g,x_{r-1}]=[c_1{\leftarrow} x_{r-1}]^{\beta_1}d_1^{\gamma_1}\cdots
  d_n^{\gamma_n}
\label{eq:gcommxrminusone}
\end{equation}
for some basic commutators $d_1<\cdots<d_n$, with $[c_1{\leftarrow}
  x_{r-1}]<d_1$. We may assume that $\gamma_i$ is positive in each
case, and less than the corresponding $N_i$ defined as
in~(\ref{eq:defofNi}). Since $[c_1{\leftarrow}x_{r-1}]$ involves at
least three generators, if any
  replacement need to be made
  they will be among the $d_i$, and none of the
  replacements nor the commutators introduced after collecting will
  be equal to $[c_1{\leftarrow} x_{r-1}]$, which has right term of
  weight one and involves three generators; thus the exponent of
  $[c_1{\leftarrow} x_{r-1}]$ in the normal form expression for
  $[g,x_{r-1}]$ is $\beta_1$. As above, this implies that
  $\beta_1\equiv 0\pmod{p^{\alpha_i}}$, and so we conclude $\beta_1=0$
  and $g\in G_{p+1}$ by induction.

\textsc{Case 4:} \textit{The commutator $c_1$ involves only the
  generators $x_{r-1}$ and~$x_r$, and
  $c_1<[x_r,x_{r-1},{}_{p-2}\,x_r]$.} We have $0\leq \beta_1 <
  p^{\alpha_{r-1}}$. We consider $[g,x_r]$; the only basic commutator
  that may need to be replaced in the expression~(\ref{eq:gcommxr}) is
  $[x_r,x_{r-1},{}_{p-1}\,x_r]$, which may appear as one of the $d_i$,
  but not as $[c_1{\leftarrow}x_r]$. If such a replacement is
  necessary, the exponent of
  $[x_r^p,x_{r-1}]$ in the normal form expression will be equal to
  $\gamma_i$, the exponent of~$d_i$ before the rewriting; this follows
  from~(\ref{eq:ypcommx}). See also \cite{struiktwo}*{Equation~(59)}. 

  If $\alpha_{r-1}<\alpha_r$, then $[g,x_r]=e$ implies that
  $\gamma_i\equiv 0 \pmod{p^{\alpha_{r-1}}}$, which contradicts our
  assumption on the $\gamma_i$ (which we assumed to be positive and
  strictly smaller than~$p^{\alpha_{r-1}}$). Thus, if
  $\alpha_{r-1}<\alpha_r$, then~(\ref{eq:gcommxr}) is already in
  normal form; since $g$ is central we must have $\beta_1\equiv 0
  \pmod{p^{\alpha_{r-1}}}$, and so we deduce $\beta_1=0$ and $g\in
  G_{p+1}$.

  If, on the other hand, $\alpha_{r-1}=\alpha_r$ then we can only
  conclude that $\gamma_i\equiv 0 \pmod{p^{\alpha_{r-1}-1}}$. Writing
  $\gamma_i = k p^{\alpha_{r-1}-1}$, then using the same argument as
  in~(\ref{eq:ytothepxraised}) we have that we will replace
  $d_i^{\gamma_i}$ with
  $[x_r,x_{r-1}]^{p\gamma_i}$ (using the fact that
  $[x_r^p,x_{r-1}]$ is of order $p^{\alpha_{r-1}-1}$).
  To write this in normal form we just need to move
  $[x_r,x_{r-1}]^{p\gamma_i}$ to the left, which introduces no
  new commutators since all other terms are already central. Thus, the
  exponent of $[c_1{\leftarrow} x_r]$ remains $\beta_1$ in the new
  expression. Again we conclude that $\beta_1\equiv 0
  \pmod{p^{\alpha_{r-1}}}$ and so $\beta_1=0$; induction now gives us
  that $g\in G_{p+1}$. 

\textsc{Case 5:} \textit{The only remaining case, $c_1 =
  [x_r,x_{r-1},{}_{p-2}\,x_r]$.}

Our assumption on~$\beta_1$ is $0\leq
  \beta_1 < p^{\alpha_{r-1}}$.
If $p=2$, then $g\equiv c_1^{\beta} = [x_r,x_{r-1}]^{\beta}
\pmod{G_3}$ with $0\leq \beta < p^{\alpha_{r-1}}$. Thus we have 
$e = [g,x_r] =
     [x_r,x_{r-1}]^{-2\beta}[x_r^2,x_{r-1}]^{\beta}$,
so by \cite{struikone}*{Theorem~4} we conclude that $-2\beta\equiv 0
\pmod{2^{\alpha_{r-1}+1}}$. From this once again we obtain that 
$\beta=0$ and $g\in G_3$. 

If $p>2$, then consider $[g,x_{r-1}]$. 
As above, the expression in~(\ref{eq:gcommxrminusone}) will be in
normal form unless one
of the commutators $d_i$ is equal to $[x_r,x_{r-1},{}_{p-1}\,x_r]$; we
now proceed as above to conclude that if $\alpha_{r-1}<\alpha_r$
then no $d_i$ needs to be replaced; and if
$\alpha_{r-1}=\alpha_r$, then we deduce that $\gamma_i\equiv 0
\pmod{p^{\alpha_{r-1}-1}}$, and so we simply replace $d_i^{\gamma_i}$
with $[x_r,x_{r-1}]^{p\gamma_i}$ and then shift this commutator to the
left, without changing the exponent of
$[c_1{\leftarrow} x_{r-1}]$. Both cases imply $\beta_1=0$
and so $g\in G_{p+1}$ by induction. 

Thus we conclude that if $g\in G_{p}\cap Z(G)$, then $g\in
G_{p+1}$. This proves that
$Z(G) = \langle x_r^{p^{\alpha_{r-1}+1}},G_{p+1}\rangle$,
as claimed. 
\end{proof}

This yields the desired result:

\begin{theorem}[cf. \cite{capable}*{Theorem~5.2}]
Let $p$ be a prime, and let $C_1,\ldots,C_r$ be cyclic $p$-groups
generated by $x_1,\ldots,x_r$, respectively; assume that the order of
$x_i$ is~$p^{\alpha_i}$ and $1\leq \alpha_1\leq\cdots\leq
\alpha_r$. If $G$ is the $p$-nilpotent product of the $C_i$,
\[ G = C_1 \amalg^{\germ N_p}\cdots\amalg^{\germ N_p} C_r,\]
then $G$ is capable if and only if $r>1$ and $\alpha_r\leq
\alpha_{r-1}+1$.
\end{theorem}

\begin{proof} Necessity follows from \cite{capable}*{Theorem
    3.12}. For sufficiency, let $K$ be the $(p+1)$-nilpotent product
      of the $C_i$, $K=C_1\amalg^{\germ N_{p+1}}\cdots\amalg^{\germ
      N_{p+1}} C_r$. By Theorem~\ref{thm:centerpplusone}, $Z(K)$ is
      generated by $x_r^{p^{\alpha_{r-1}+1}}$ and $K_{p+1}$. Since
      $\alpha_r\leq \alpha_{r-1}+1$, the former is trivial, so
      $Z(K)=K_{p+1}$. Thus $K/Z(K) = K/K_{p+1}\cong G$, so $G$ is
      capable.
\end{proof}


\section*{Reference}

\begin{biblist}

\bib{hurley}{article}{
   author = {Hurley, T. C.},
    title = {Identifications in free groups},
     date = {1987},
  journal = {J. Pure Appl. Algebra},
   volume = {48},
   number = {3},
    pages = {249\ndash 261},
   review = {\MR{89a:20025}},
}

\bib{capable}{article}{
    author={Magidin, Arturo},
     title={Capability of nilpotent products of cyclic groups},
   journal={J. Group Theory},
    volume={8},
    number={4},
      year={2005},
   pages={431\ndash 452},
  review={\MR{2006c:20073}},
}

\bib{struikone}{article}{
    author={Struik, Ruth~Rebekka},
     title={On nilpotent products of cyclic groups},
      date={1960},
   journal={Canad. J. Math.},
    volume={12},
     pages={447\ndash 462},
    review={\MR{22:\#11028}},
}

\bib{struiktwo}{article}{
    author={Struik, Ruth~Rebekka},
     title={On nilpotent products of cyclic groups II},
      date={1961},
   journal={Canad. J. Math.},
    volume={13},
     pages={557\ndash 568},
    review={\MR{26:\#2486}},
}
\bib{naturalorder}{article}{
   author={Waldinger, Hermann V.},
    title={A natural linear ordering of basic commutators},
     date={1961},
  journal={Proc. Amer. Math. Soc.},
   volume={12},
    pages={140\ndash 147},
   review={\MR{26:\#2485}},
}

\bib{witts}{article}{
   author={Waldinger, Hermann V.},
    title={On extending {W}itt's formula},
     date={1967},
  journal={J. Algebra},
   volume={5},
    pages={41\ndash 58},
   review={\MR{34:\#228}},
}



\bib{basiccomms}{article}{
   author = {Ward, M. A.},
    title = {Basic commutators},
   journal = {Philos. Trans. Roy. Soc. London Ser. A},
     date = {1969},
   volume = {264},
    pages = {343\ndash 412},
   review = {\MR{40:\#4379}},
}
\end{biblist}
\end{document}